\documentclass[12pt]{article}
\usepackage{setspace}
\usepackage{amssymb}
\usepackage{amsthm}
\usepackage{amsfonts}
\usepackage{natbib}
\bibpunct{(}{)}{;}{a}{}{,}
\begin{document}

\title{Moments of Sums of Independent and Identically Distributed Random Variables}

\author{\\ \\ \\ Daniel M. Packwood \\ \\ Department of Chemistry \\ Graduate School of Science\\ Kyoto University, Kyoto, 606-8502, Japan}
\maketitle

\newpage

\begin{abstract}

We present an analytic method for computing the moments of a sum of independent and identically distributed random variables. The limiting behavior of these sums is very important to statistical theory, and the moment expressions that we derive allow for it to be studied relatively easily. We show this by presenting a new proof of the central limit theorem and several other convergence results. 


\end{abstract}

\newpage

\section{Introduction}

Sums of independent and identically distributed random variables often appear in statistics and probability. Their major role in statistics is to provide a mathematical model for the estimation of a population mean from a sample. We write such sums as

\[
S_n = \sum_{i=1}^{n} X_i,
\]

\noindent
where $X_1, X_2,\ldots,X_n$ are independent and identically distributed (\textit{iid}) random variables distributed as a random variable $X$, where $E(X) = \mu$ and $\mbox{var}(X) = \sigma^2$. The limiting behavior of these sums is very significant. For example, the central limit theorem says that under certain conditions

\[
\frac{S_n - n\mu}{\sigma \sqrt{n}} \Rightarrow Z
\]

\noindent
as $n \rightarrow \infty$, where $Z \sim N(0,1)$ and $\Rightarrow$ denotes weak convergence (convergence in distribution). Because of this kind of limiting behavior, sums of \textit{iid} random variables form the basis of asymptotically normal statistical estimation and hypothesis testing.

Probabilists have been interested in the moments of sums of random variables since the early part of last century. Khinchine's 1923 paper appears to make the first significant contribution to this problem [1]. It provides inequalities for the moments of a sum of Bernoulli random variables. In 1970, Rosenthal generalised Khinchine's result to the case of positive or mean-zero random variables [2]. Further refinements to these bounds have been made by Latala and Hitczenko, Montogomery-Smith and Oleszkiewiez in more recent times [3, 4, 5]. Nonetheless, there appears to be no equalities available for the case of a sum of \textit{iid} random variables. Amongst other uses, expressions for the moments should allow for the limiting behavior of $S_n$ to be studied relatively easily.

This article presents an analytic method for computing each moment of a sum of \textit{iid} random variables. As an application, we will take limits of these expressions and obtain the central limit theorem and some other results. These results are probabily easy to obtain by other methods, however they do not appear to have been highlighted in the literature. We  prove a limit theorem for the case where $X$ has an \emph{asymmetric distribution} about 0 (i.e., $P(-\infty < X < - c) \neq P(c < X < \infty)$), namely

\[
\frac{S_n}{n} \Rightarrow H_\mu,
\] 

\noindent
where $H_\mu$ is a random variable with a unit-step distribution,

\[
P(H_\mu < h) = \left\{\begin{array}{ll}
0 & h < \mu \\
1 & h \geq \mu \\
\end{array} \right.
\] 

\noindent
We also obtain a generalisation of the weak law of large numbers,

\[
\lim_{n \rightarrow \infty}P(|S_n^p/n^p - \mu^p| > \epsilon) = 0, 
\]

\noindent
for all $\epsilon > 0$ and $p = 1, 2, \ldots$, and a generalisation of the strong law of large numbers for when $X$ has a symmetric distribution about $\mu = 0$,

\[
S_n^p/n^p \rightarrow 0,
\]

\noindent
where $\rightarrow$ indicates almost sure convergence. The conditions under which these results hold are discussed in the main text.

Section 2 provides a formula for the moments and its proof, and also shows how the formula is applied. Section 3 studies the limits of these formulas section 4 presents some final remarks.

\section{Formula for the moments}

We present the formula for the moments as follows. Fix a $p \in \{1,2,\ldots\}$ and define

\[
\mathcal{Q}^p = \left\{\begin{array}{ll}
E(X^r)E(X^s)\cdots E(X^t) : & r, s, \ldots, t \in \{1,2,\ldots, p\} \\  
														& r + s + \cdots + t = p
										\end{array} \right\}.
\]

\noindent
The $p$th moment of $S_n = \sum_{i=1}^n X_i$ is

\begin{equation}
E(S_n^p) = \sum_{q_i \in \mathcal{Q}^p} a_i q_i,
\end{equation}

\noindent
where, for $q_i = E(X^{p_1})E(X^{p_2})\cdots E(X^{p_m})$,

\begin{equation}
a_i = \frac{1}{l_1!l_2!\cdots l_h!}\frac{n!}{(n - m)!}\frac{p!}{p_1! p_2! \cdots p_m!}.
\end{equation}

\noindent
In Equation (1), $h$ is the number of distinct constants in the sequence \allowbreak$\{p_1, p_2, \ldots, p_m\}$, $l_1$ the number of elements equal to the first constant, $l_2$ the number equal to the second constant, $\ldots$, and $l_h$ the number equal to the $h$th constant (e.g., for the sequence $\{1,1,1,1\}$, $h = 1$ and $l_1 = 4$, and for $\{1,2,2,2\}$, $h = 2$, $l_1 = 1$, and $l_2 = 3$). Note that $E(S_n^p) < \infty$ if an only if $E(X^\alpha) < \infty$ for all positive integers $\alpha \leq p$.

For example, the third moment works out to be

\begin{eqnarray}
E(S_n^3) &=& \left(\frac{1}{1!}\frac{n!}{(n-1)!}\frac{3!}{3!}\right)E(X^3) + \left(\frac{1}{1!1!}\frac{n!}{(n-2)!}\frac{3!}{2!1!}\right)E(X^2)E(X) \nonumber \\
				 &+& \left(\frac{1}{3!}\frac{n!}{(n-3)!}\frac{3!}{1!1!1!}\right) E(X)^3 \nonumber \\
				 &=& nE(X^3) + 3n(n-1)E(X^2)E(X) + n(n-1)(n-2)E(X)^3 \nonumber
\end{eqnarray}

The proof of Equations (1) and (2) is a lengthly application of various combinatorial formulas. 

\begin{proof}

The $p$th moment of $S_n$ is

\begin{equation}
E(S_n^p) = \underbrace{\sum_{i=1}^n \sum_{j=1}^n \cdots \sum_{k=1}^n}_{p \mbox{ sums}} E(X_i X_j \cdots X_k),
\end{equation}
Because $X_1,\ldots,X_n$ are \textit{iid} random variables, each term in the sum can be factored into the form 

\[
E(X^{p_1}) E(X^{p_2}) \cdots E(X^{p_n}), 
\]

\noindent
where $p_1, p_2, \ldots, p_n$ are positive integers that sum to $p$. Equation (3) can therefore be written as Equation (1). To determine the constants $a_1, a_2, \ldots$ fix a positive value of $p$ and positive non-zero integers $p_1, p_2,\ldots, p_m$ such that $p_1 + p_2 + \ldots + p_m = p$. For each choice of $i,j,\ldots,k \in \{1,2,\ldots,n\}$, with $i  \neq j \neq, \ldots \neq k$, define the collection

\[
\mathcal{C}_{i_{p_1},j_{p_2},\ldots,k_{p_m}} = \left\{\begin{array}{cc}
	\mbox{all unique permutations of the sequence } \\
	\underbrace{X_i, X_i, \ldots X_i}_{p_1 \mbox{ objects}}, \underbrace{X_j, X_j, \ldots X_j}_{p_2 \mbox{ objects}} 		,\ldots, \underbrace{X_k, X_k, \ldots X_k}_{p_m \mbox{ objects}} \\
	\end{array}
	\right\}	
\]

\noindent
Each of the sequences in $\mathcal{C}_{i_{p_1},j_{p_2},\ldots,k_{p_m}}$ corresponds to an expectation of the form $E(X_i^{p_1} X_j^{p_2} \cdots X_k^{p_m})$. These are each equivalent to

\[
q_m = E(X^{p_1})E(X^{p_2})E(X^{p_m})
\] 

\noindent
The subscript $m$ on $q_m$ is only for bookkeeping purposes. The number of elements in $\mathcal{C}_{i_{p_1},j_{p_2},\ldots,k_{p_m}}$ is the number of ways in which $q_m$ can be created from $p_1$ copies of $X_i$, $p_2$ copies of $X_j$, $\ldots$, $p_m$ copies of $X_k$. For any choice of $i,j,\ldots,k$, this number is

\begin{eqnarray}
\vert \mathcal{C}_{i_{p_1},j_{p_2},\ldots,k_{p_m}} \vert &=& \frac{p!}{p_1! p_2! \cdots p_m!}.
\end{eqnarray}

\noindent
The number of elements in the union

\begin{equation}
\mathcal{Q}_{p_1,\ldots,p_m}^p = \bigcup_{\begin{array}{c}
			i,j,\ldots,k \in \{1,\ldots, n\} \\
			i \neq j \neq \cdots \neq k 
											 \end{array}}
			\mathcal{C}_{i_{p_1},j_{p_2},\ldots,k_{p_m}}			
\end{equation}

\noindent
is the number of ways $q_m$ can be created from the $p_1$ copies of $X_i$, $p_2$ copies of $X_j$, $\ldots$, $p_m$ copies of $X_k$, $n$ choices of $i$, $n-1$ choices of $j$, $\ldots$, $n-(m-1)$ choices of $k$. This gives the identity

\begin{equation}
a_m = \vert \mathcal{Q}_{p_1,\ldots,p_m}^p \vert.
\end{equation} 

\noindent
If there are no equalities between the constants $p_1, p_2, \ldots, p_m$, then the collections $\{\mathcal{C}_{i_{p_1},j_{p_2},\ldots,k_{p_m}}\}$ are mutually exclusive and so $\vert \mathcal{Q}_{p_1,\ldots,p_m}^p \vert$ is given by Equation (4) times the number of collections in the union, namely

\[
n(n-1)\cdots (n-(m-1)) \frac{p!}{p_1! p_2! \cdots p_m!}.
\]

\noindent
If there are equalities between $p_1, p_2, \ldots p_m$, then some of the collections in $\{\mathcal{C}_{i_{p_1},j_{p_2},\ldots,k_{p_m}}\}$ are equivalent. For example, if $p_1 = p_2 = p_3$, then

\begin{eqnarray*}
\mathcal{C}_{a_{p_1},b_{p_2},c_{p_3},\ldots,k_{p_m}} 
		&=& \mathcal{C}_{b_{p_1},c_{p_2},a_{p_3},\ldots,k_{p_m}} \nonumber \\
		&=& \mathcal{C}_{c_{p_1},a_{p_2},b_{p_3},\ldots,k_{p_m}} \nonumber \\
		&=& \mathcal{C}_{c_{p_1},b_{p_2},a_{p_3},\ldots,k_{p_m}} \nonumber \\
		&=& \mathcal{C}_{a_{p_1},c_{p_2},b_{p_3},\ldots,k_{p_m}} \nonumber \\
		&=& \mathcal{C}_{b_{p_1},a_{p_2},c_{p_3},\ldots,k_{p_m}} \nonumber \\
\end{eqnarray*} 

\noindent
In the above equalities, the indices given by the ellipses and $k_{p_m}$ do not change. Furthermore, if we also had $p_4 = p_5 = p_6 = p_7$ in addition to $p_1 = p_2 = p_3$, then we could also write an additional $4!$ equivalent collections for each of the $3!$ equivalent collections given above, leading to $4!3!$ equivalent collections in the union in Equation (5). In general, if there are $h$ distinct constants in the sequence $\{p_1, p_2, \ldots p_m\}$, with $l_1$ of the elements equal to one of these constants, $l_2$ equal to another of these constants, \ldots, $l_h$ equal to the remaining of these constants, then $l_1! l_2! \cdots l_h!$ of the collections in the union in Equation (5) are equivalent. The number of elements in $\mathcal{Q}_{p_1,\ldots,p_m}^p$ for the general case is therefore found by taking the mutually exclusive result given above and dividing through by $l_1! l_2! \cdots l_h!$, namely

\[
\vert \mathcal{Q}_{p_1,\ldots,p_m}^p \vert = \frac{1}{l_1! l_2! \cdots l_h!} \frac{n!}{(n-m)!} \frac{p!}{p_1! p_2! \cdots p_m!}.
\]

\noindent
Equation (6) then gives Equations (1) and (2).

\end{proof}

\section{Application: Limit theorems for $S_n$}

The moment equations (1) and (2) are convenient for investigating the behavior of $S_n$ as $n \rightarrow \infty$. We can do this with the \emph{method of moments}: Let $F_n$ be the distribution function of $S_n$ and $F$ a distribution function with moments $m_1,$ $m_2$, $\ldots$. Assume that $F$ is determined completely by its moments. If $E(S_n^k) \rightarrow m_k$ for every $k \geq 0$, then $F_n \rightarrow F$ (see e.g., [6]). 

\subsection{The classical central limit theorem}

The classical central limit theorem is particularly easy to derive from the moment equations. Without loss of generality we will assume that $X$ has a symmetric distribution about 0. Equations (1) and (2) simplify significantly in this case because all terms involving odd moments of $X$ are identically zero. Moreover, all odd moments of $S_n$ are also identically zero. This can be seen from Equation (3). For odd $p$, at least one of $p_1$, $p_2$, $\ldots$, $p_n$ will be odd, and so each term in Equation (3) will contain at least one odd moment of $X$ and will be zero. The second, fourth and sixth moments for the symmetric case work out to be

\begin{equation}
E(S_n^2) = n\sigma^2
\end{equation}

\begin{equation}
E(S_n^4) = n\mu_4 + 3n(n-1)\sigma^4
\end{equation}

\begin{equation}
E(S_n^6) = n\mu_6 + 15n(n-1)\mu_4\sigma^2 + 15n(n-1)(n-2)\sigma^6
\end{equation}

\noindent
where $\mu_k = E(X^k)$ and $\sigma^2 = E(X^2)$. After some manipulations, the moments of $S_n/(\sigma\sqrt{n})$ take on the form

\begin{equation}
E\left(\frac{S_n^p}{n^{p/2}\sigma^p}\right) = (p - 1)!! + \beta_n(p),
\end{equation}

\noindent
where $k!! = k(k-2)(k-4)\cdots1$ is the double factorial and $(p-1)!!$ the $p$th moment of a standard normal random variable. $\beta_n(p)$ can be regarded as a non-normal `correction term' for the $p$th moment. While there appears to be no general formula for the correction terms, the first few are

\begin{equation}
\beta_n(2) = 0
\end{equation}

\begin{equation}
\beta_n(4) = \frac{1}{n}\left(\frac{\mu_4}{\sigma^4} - 3\right)
\end{equation}

\begin{equation}
\beta_n(6) = \frac{1}{n}\left(\frac{15\mu_4}{\sigma^4} - 45\right) + \frac{1}{n^2}\left(\frac{\mu_6}{\sigma^6} - \frac{15\mu_4}{\sigma^4} + 30\right)
\end{equation}

\noindent
In general, $\beta_n(p) \rightarrow 0$ for all $p \geq 0$, and so $E(S_n^p/(n^{p/2}\sigma^p)) \rightarrow (p - 1)!!$ completely. The classical central limit theorem then follows from the fact that the normal distribution is completely determined by its moments.

The above assumes that all moments of the random variable $X$ are finite. While this covers many important cases, other proofs of the central limit theorem do not make this assumption and are less restrictive.

\subsection{Convergence of a sum of asymmetric random variables}

Now we consider the case where $X$ has an asymmetric distribution about 0. Under this assumption, the odd moments of $S_n$ are no longer zero and there are additional terms in the expressions for the even moments. The first three moments are

\begin{equation}
E(S_n) = n\mu
\end{equation}

\begin{equation}
E(S_n^2) = n\mu_2 + n(n-1)\mu^2
\end{equation}

\begin{equation}
E(S_n^3) = n\mu_3 + 3n(n-1)\mu_2\mu + n(n-1)(n-2)\mu^3,
\end{equation}

\noindent
and after some manipulations the moments of $S_n/n$ work out to be

\begin{equation}
E\left(\frac{S_n^p}{n^p}\right) = \mu^p + \alpha_n(p)
\end{equation}

\noindent
where $\alpha_n(p)$ is another correction term. The first three correction terms are

\[
\alpha_n(1) = 0
\]

\[
\alpha_n(2) = \frac{1}{n}\left(\mu_2 - \mu^2\right)
\]

\[
\alpha_n(3) = \frac{1}{n^2}\left(\mu_3 - 3\mu_2\mu + 2\mu^3\right) + \frac{1}{n}\left(3\mu_2\mu - 3\mu^3\right)
\]

\noindent
In the general case, we have $\alpha_n(p) \rightarrow 0$ for all $p \geq 0$, and hence

\begin{equation}
E\left(\frac{S_n^p}{n^p}\right) \rightarrow \mu^p
\end{equation}

\noindent
The sequence of moments $\mu$, $\mu^2$, $\mu^3$, $\ldots$ determine a probability distribution function with a unit-step at $\mu$, namely

\[
P(H_\mu < h) = \left\{\begin{array}{ll}
0 & h < \mu \\
1 & h \geq \mu \\
\end{array} \right.
\]

\noindent
This can be checked directly by expanding the characteristic function of $H_\mu$ into its Taylor's series and making the substitution $E(H_\mu^k) = \mu^k$.

We have therefore proven the theorem mentioned in the introduction, namely

\begin{equation}
\frac{S_n}{n} \Rightarrow H_\mu.
\end{equation}

\noindent
Our proof has again required all moments of $X$ to be finite. Loosely speaking, Equation (19) says that the probability density of $S_n/n$ concentrates into a delta-like point mass as $n \rightarrow \infty$. This kind of behavior is apparent from other asymptotic results as well. For example, it is well-known (and easy to show from Equations (14) - (16)) that for very large $n$ $S_n/n$ is approximately $N(\mu,\sigma^2/n)$. This result suggests that as $n \rightarrow \infty$, the probability density of $S_n/n$ `converges to a delta function' centered at $\mu$. The novelty of Equation (19) is that it puts this concept on precise mathematical footing.

\subsection{Generalisations of the laws of large numbers}

Equation (18) in the previous section says that for asymmetric $X$ the $p$th moment of $S_n/n$ converges to $\mu^p$. This is also true for symmetric $X$, in which $\mu = 0$. This can be seen from Equations (7) - (9), which show that the highest power of $n$ that appears in the numerator of $E(S_n^p)$ is $p/2$, and so for large $n$

\begin{equation}
E\left(S_n^p/n^p\right) = O\left(n^{-p/2}\right) \rightarrow 0.
\end{equation}

\noindent
We will now prove a generalisation of the weak law of large numbers, namely

\begin{equation}
\lim_{n \rightarrow \infty} P\left(\left|S_n^p/n^p - \mu^p\right| > \epsilon\right) = 0
\end{equation}

\noindent
for all $\epsilon > 0$. The usual weak law is the case $p = 1$. As with the usual law, the proof is based on the Markov inequality:

\begin{equation}
P\left(\left|S_n^p/n^p - \mu^p\right| > \epsilon\right) \leq \frac{1}{\epsilon^2}E((S_n^p/n^p - \mu^p)^2).
\end{equation}

\noindent
For the case of symmetric $X$, $\mu = 0$ and so by Equation (20) the right-hand side of Equation (22) goes to zero. For the asymmetric case,

\[
E((S_n^p/n^p - \mu^p)^2) = E(S_n^{2p}/n^{2p}) + \mu^{2p} - 2\mu^p E(S_n^p/n^p) \rightarrow 0
\]

\noindent
by Equation (18). Taking the limit of Equation (22) then completes the proof.

A generalisation of the strong law of large numbers can be proven for the case where $X$ is symmetric. From Equations (7) - (9) we can see that

\[
E\left(S_n^{2p}\right) = \sum_{m=1}^p c_m n^m \leq Cn^p,
\]

\noindent
where $c_1$, $c_2$, $\ldots$, $c_p$ are constants and $C = \max_{m = 1,2,\ldots, p} \left\{c_m\right\}$. Summing Equation (22) over $n$ therefore gives

\begin{eqnarray}
\sum_{n=1}^{\infty} P\left(\left|S_n^p/n^p\right| > \epsilon \right) &\leq& \frac{1}{\epsilon^2}\sum_{n=1}^{\infty} E(S_n^{2p}/n^{2p}) \nonumber \\
&\leq& \frac{C}{\epsilon^2} \sum_{n=1}^{\infty} n^{-p} < \infty \nonumber
\end{eqnarray}

\noindent
By the Borel-Cantelli lemma,

\begin{equation}
P\left(\omega \mbox{ that are in infinitely many } \left\{\left|S_n^p/n^p\right| > \epsilon\right\}\right) = 0,
\end{equation}

\noindent
and so there exists an $n_0$ such that for $n > n_0$, $P\left(\left|S_n^p/n^p\right| < \epsilon \right) = 1$. Letting $\epsilon \rightarrow 0$, we find that

\begin{equation}
S_n^p/n^p \rightarrow 0,
\end{equation}

\noindent
where $\rightarrow$ indicates \textit{almost sure} convergence. It does not appear possible to apply this method to the case of asymmetric $X$.

\section{Final remarks}

We have presented an analytic method of computing the moments of a sum of \textit{iid} random variables and used it to derive the central limit theorem and several other new limiting results. The word \emph{`method'} should be emphasised because Equations (1) and (2) only give a systematic way of calculating the moments, rather than general formula for the moments explicitly. Indeed, from Equations (7) to (9) and Equations (14) to (16), the moments do not appear to possess any general structure. The method looks useful for numerical computation of higher moments, although we have not explored this possiblity here.

Outside of statistics, the major application of a sum of \textit{iid} random variables is in the study of random walks, which play a fundamental role in stochastic process theory and statistical physics (see e.g., [7, 8]). We have discussed continuous time random walks in a physical context elsewhere [9, 10, 11].

\vspace{10 mm}
\noindent
\textbf{Acknowledgements} \\
D.M.P is supported by a Japan Society for the Promotion of Science Postdoctoral Fellowship. An anonymous reviewer of an earlier version of this manuscript is thanked for helpful suggestions.

\vspace{10 mm}
\noindent
\textbf{References} \\

\noindent
[1] Khintchine, A. Uber dyadische bruche. \textit{Mathematische Zeitschrift}, \textbf{18}, 1923, 109 - 116. \\

\noindent
[2]	Rosenthal, H. P. On the subspaces of $L^p$ ($p > 2$) spanned by sequences of independent random variables. \textit{Israel Journal of Mathematics}, \textbf{8}, 1970, 273 - 303. \\

\noindent
[3] Latala, R. Estimation of moments of sums of independent real random variables. \textit{The Annals of Probability}, \textbf{25}, 1997, 1502 - 1513. \\

\noindent
[4] Hitczenko, P., Montgomery-Smith, S. J., Oleszkiewiez, K. Moment inequalities for sums of certain independent symmetric random variables. \textit{Studia Mathematicia}, \textbf{123}, 1997, 15 - 45. \\

\noindent
[5] Hitczenko, P., Montgomery-Smith, S. J. Measuring the magnitude of sums of independent random variables. \textit{The Annals of Probability}, \textbf{29}, 2001, 447 - 466. \\

\noindent
[6] Galambos, J. Advanced Probability Theory. 1995, Marcel Dekker, New York, N.Y, U.S.A. \\

\noindent
[7] Durrett, R. Probability: Theory and Examples. 2010, Cambridge University Press, New York, N.Y, U.S.A. \\

\noindent
[8] van Kampen, N. G. Stochastic Processes in Physics and Chemistry. 1987, Elsevier, New York, N.Y., U.S.A. \\

\noindent
[9] Packwood, D. M., Tanimura, Y. Non-Gaussian stochastic dynamics of spins and oscillators: A continuous-time random walk approach. \textit{Physical Review E}, \textbf{84}, 2011, 61111 - 61124. \\

\noindent
[10] Packwood, D. M., Tanimura, Y. Phase relaxation functions for a continuous time random walk (working title). In preparation. \\

\noindent
[11] Packwood, D. M. Relaxation function for a continuous time random walk (Japanese). \textit{Proceedings of the 8th Mathematics Conference for Young Researchers}, 2012, Hokkaido University. In press. \\

\end{document}